# PARAEXPONENTIALS, MUCKENHOUPT WEIGHTS, AND RESOLVENTS OF PARAPRODUCTS

CRISTINA PEREYRA AND LESLEY WARD

ABSTRACT. We analyze the stability of Muckenhoupt's $\mathbf{RH}_p^d$ and $\mathbf{A}_p^d$ classes of weights under a nonlinear operation, the $\lambda$-operation. We prove that the dyadic doubling reverse Hölder classes $\mathbf{RH}_p^d$ are not preserved under the $\lambda$-operation, but the dyadic doubling $A_p$ classes $\mathbf{A}_p^d$ are preserved for $0 < \lambda < 1$. We give an application to the structure of resolvent sets of dyadic paraproduct operators.

## 1. INTRODUCTION

The Muckenhoupt classes of weights consist of positive locally integrable functions satisfying certain integrability conditions on intervals; see Section 2. They arise, for instance, in connection with the boundedness of the Hilbert transform and maximal function operators on $L^p$ spaces.

Weights in the Muckenhoupt classes are often represented as exponentials of functions in **BMO**; see [GC-RF]. For weights $\omega(x)$ defined on an interval $J$, with mean value one on $J$, a different representation, sometimes called the *paraexponential*, was introduced in [FKP]. The correspondence is realized by an infinite product, namely:

$$(1.1) \qquad \omega(x) = \prod_{I \in \mathcal{D}(J)} \left(1 + b_I h_I(x)\right), \qquad b_I = \frac{\langle \omega, h_I \rangle}{m_I \omega};$$

where $\mathcal{D}(J)$ denotes the intervals in the dyadic decomposition of the interval $J$, $\{h_I\}_{I \in \mathcal{D}(J)}$ are Haar functions, $\langle \cdot, \cdot \rangle$ is the scalar product in $L^2$, and $m_I \omega$ is the mean value of $\omega$ on $I$.

It was proved in [FKP] that if there is an $\varepsilon > 0$ such that $|b_I h_I| < 1 - \varepsilon$ for all $I \in \mathcal{D}(J)$, then the partial products converge to a weight in dyadic doubling $\mathbf{A}_\infty^d$ if and only if $b = \sum_{I \in \mathcal{D}(J)} b_I h_I$ is a function of dyadic bounded mean oscillation, $\mathbf{BMO}^d$.

1991 *Mathematics Subject Classification.* 42c, 47b; Secondary 41.
*Key words and phrases.* Muckenhoupt weights, reverse Hölder $RH_p$, $A_p$, doubling weights, dyadic paraproducts, paraexponentials.
Research supported in part by (CP) NSF grant #DMS-93-04580 and (LW) at MSRI by NSF grant #DMS-90-22140

A larger dictionary relating properties of $\omega$ and $b$ can be found in [B]. In particular the dyadic doubling Muckenhoupt $\mathbf{RH}_p^d$ and $\mathbf{A}_p^d$ classes can be characterized by summation conditions on $b$.

We will consider only weights $\omega$ defined on an interval $J$, with mean value one on $J$. For such a weight we have, at least formally, the product representation (1.1). We are interested in the effect of multiplying the coefficient $b_I$ in each factor of the product by the same number $\lambda$. Define

$$(1.2) \qquad \omega_\lambda(x) = \prod_{I \in \mathcal{D}(J)} \left(1 + \lambda b_I h_I(x)\right).$$

We call the mapping that sends $\omega \mapsto \omega_\lambda$ the $\lambda$-*operation*.

To guarantee the convergence of the products we will only consider the case $-1 \leq \lambda \leq 1$. If $\omega$ is a weight in $\mathbf{A}_\infty^d$ and $-1 \leq \lambda \leq 1$, then $\omega_\lambda$ is also in $\mathbf{A}_\infty^d$; see [FKP].

The $\lambda$-operation seems, at first sight, very similar to taking a $\lambda$ power of the weight, $\omega^\lambda$. The first two terms in the Taylor expansion of $\omega^\lambda = e^{\lambda b}$ coincide with the first two terms obtained by expanding the infinite product in (1.2), but the third terms differ. Both the $\mathbf{RH}_p^d$ and the $\mathbf{A}_p^d$ classes are preserved under taking powers $\omega \mapsto \omega^\lambda$, for $0 \leq \lambda \leq 1$.

**Question:** Are the $\mathbf{RH}_p^d$ and $\mathbf{A}_p^d$ classes preserved under the $\lambda$-operation?

The answer is negative for the $\mathbf{RH}_p^d$ classes, and positive for the $\mathbf{A}_p^d$ classes.

**Theorem 1.1.** *For each $p > 1$, there exist a weight $\omega$ in $\mathbf{RH}_p^d$ and a number $\lambda \in (0,1)$ such that $\omega_\lambda$ is not in $\mathbf{RH}_p^d$.*

For $p > (1 - \log 2)^{-1} > 1$, examples with positive $\lambda$'s are given in [P2]. In this paper we give examples with positive $\lambda$'s for all $p > 1$. Examples for all $p > 1$, with $-1 \leq \lambda < 0$, are given in [P1].

A related result in [P2] is that if $\omega \in \mathbf{RH}_p^d$, then there is a $q > 1$ such that $\omega_\lambda \in \mathbf{RH}_q^d$ for all $-1 \leq \lambda \leq 1$.

**Theorem 1.2.** *Given a weight $\omega \in \mathbf{A}_p^d$, then $\omega_\lambda \in \mathbf{A}_p^d$ for all $1 \leq p \leq \infty$ and all $0 \leq \lambda \leq 1$.*

However, the $\mathbf{A}_p^d$ spaces for $1 \leq p < \infty$ are not preserved by the $\lambda$-operation for negative $\lambda$; this can be shown by examples like those in [P1].

The $\lambda$-operation appears naturally in the study of the resolvents of the dyadic paraproduct [P2]. The dyadic paraproduct is a bilinear operator that appears in different guises in harmonic analysis, often replacing the ordinary product [M, Ch, D]. The dyadic paraproduct $\pi_b$ associated to a function $b$ is defined by

$$(1.3) \qquad \pi_b f(x) = \sum_{I \in \mathcal{D}(J)} m_I f \, b_I h_I(x),$$

where $m_I f$ is the mean value of $f$ on the interval $I$, and $b_I = \langle b, h_I \rangle$. The paraproduct $\pi_b$ is bounded on $L^p(J)$ if and only if $b \in \mathbf{BMO}^d$. The conditions on $b$ that guarantee the existence of a bounded inverse of $(I - \lambda \pi_b)$ on $L^p_o(J) = \{f \in L^p(J) : \int_J f = 0\}$ are described in terms of properties of the weight $\omega_\lambda$. For doubling weights, the conditions reduce to $\omega_\lambda \in \mathbf{RH}^d_p$; see [P1]. More precisely:

**Theorem 1.3 (P).** *Let $b$ be a function in $\mathbf{BMO}^d$, take $\varepsilon > 0$, and let $\lambda$ be a real number such that $|\lambda b_I h_I| \leq (1 - \varepsilon)$ for all $I \in \mathcal{D}(J)$, where $b_I = \langle b, h_I \rangle$. Then $(I - \lambda \pi_b)^{-1}$ exists and is bounded on $L^p_o(J)$ if and only if the weight $\omega_\lambda(x) = \prod_{I \in \mathcal{D}(J)} \left(1 + \lambda b_I h_I(x)\right)$ is in $\mathbf{RH}^d_p$.*

**Question:** If $(I - \lambda \pi_b)^{-1}$ exists and is bounded on $L^p$ for $\lambda = 1$, is the same true for all $\lambda \in [0, 1]$?

Since the resolvent set of $\pi_b$ is open, there are neighborhoods of $\lambda = 0$ and $\lambda = 1$ for which this holds. By the discussion above, the question is equivalent to asking whether the $\mathbf{RH}^d_p$ classes are preserved under the $\lambda$-operation. Theorem 1.1 shows that the answer is negative. Rephrasing the theorem we have:

**Theorem 1.4.** *For each $p > 1$, there exist a number $\lambda \in (0, 1)$, a function $b \in \mathbf{BMO}^d$, and a number $\varepsilon > 0$ such that $|b_I h_I| \leq 1 - \varepsilon$ for all $I \in \mathcal{D}(J)$, the operator $(I - \pi_b)^{-1}$ exists and is a bounded operator on $L^p_o(J)$, but the operator $(I - \lambda \pi_b)$ is not invertible as an operator from $L^p_o(J)$ into $L^p_o(J)$.*

In Section 2 we give notation, definitions, and some preparatory lemmas. In Section 3 we prove Theorem 1.2, and in Section 4 we construct the examples which establish Theorem 1.1.

As customary, $C$ denotes a constant that may change from line to line.

## 2. Preliminaries

In this section we introduce the dyadic intervals and the Haar basis; and we define dyadic $\mathbf{BMO}^d$ in terms of a Carleson condition. Next we define the dyadic doubling Muckenhoupt classes $\mathbf{A}^d_\infty$, $\mathbf{RH}^d_p$, and $\mathbf{A}^d_p$; and we recall some of their properties. We state the Fefferman-Kenig-Pipher Product Representation Theorem for $\mathbf{A}^d_\infty$ weights and Buckley's Theorem characterizing $\mathbf{A}^d_p$ weights via summation conditions. Finally we define the $\lambda$-operation, and give a convexity lemma used in the proof of Theorem 1.2.

Let $\mathcal{D}$ denote the family of all dyadic subintervals of $[0, 1]$, in other words all intervals of the form $(j 2^{-k}, (j+1) 2^{-k}]$, $j, k$ integers, $0 \leq k$, $0 \leq j \leq 2^k - 1$. Given any interval $J$, $\mathcal{D}(J)$ denotes the family of dyadic subintervals of $J$. Given an interval $J$ we denote its left and right halves respectively by $J_l$ and $J_r$. An interval $\widetilde{I}$ is the *parent* of an interval $I$ if $I$ is $\widetilde{I}_l$ or $\widetilde{I}_r$.

The *Haar function* associated to an interval $I$ is given by $h_I(x) = |I|^{-1/2}\big(\chi_{I_r}(x) - \chi_{I_l}(x)\big)$; here $\chi_I$ denotes the characteristic function of the interval $I$. The set of Haar functions indexed by $\mathcal{D}(J)$ forms a basis of $L_o^2(J) = \{f \in L^2(J) : \int_J f = 0\}$; see [H].

A locally integrable function $b$ on $[0,1]$ is in the space of *dyadic bounded mean oscillation* $\mathbf{BMO}^d$ if there is a constant $C$ such that $\int_J |b(x) - m_J b|^2\, dx \leq C\,|J|$ for all $J \in \mathcal{D}$, where $m_J b = \frac{1}{|J|}\int_J b$. The function $b$ is in $\mathbf{BMO}^d$ if and only if the *Carleson condition* on the Haar coefficients $b_I = \langle b, h_I\rangle$ of $b$ holds: $\sum_{I \in \mathcal{D}(J)} b_I^2 \leq C\,|J|$ for all $J \in \mathcal{D}$, with a constant $C$ independent of $J$; see [Ch, M].

## 2.1. Dyadic weights.

We consider weights defined on the interval $J_0 = [0,1]$.

A *dyadic doubling weight* $\omega$ is a positive locally integrable function such that $\int_{\widetilde{I}} \omega \leq C \int_I \omega$ for all intervals $I \in \mathcal{D}$, where $\widetilde{I}$ is the parent of $I$ and $C$ is a constant independent of $I$.

A weight $\omega$ is in the *dyadic doubling $A_\infty$ class* $\mathbf{A}_\infty^d$ if $\omega$ is dyadic doubling, and there is a constant $C$ such that $\frac{1}{|I|}\int_I \omega \leq C \exp\left(\frac{1}{|I|}\int_I \log \omega\right)$, for all $I \in \mathcal{D}$.

A weight $\omega$ is in the *dyadic doubling reverse Hölder class* $\mathbf{RH}_p^d$, for $1 < p < \infty$, if $\omega$ is dyadic doubling, and there is a constant $C$ such that $\left(\frac{1}{|I|}\int_I \omega^p\right)^{1/p} \leq C\,\frac{1}{|I|}\int_I \omega$, for all $I \in \mathcal{D}$.

A weight $\omega$ is in the *dyadic doubling $A_p$ class* $\mathbf{A}_p^d$, for $1 < p < \infty$, if $\omega$ is dyadic doubling, and there is a constant $C$ such that $\left(\frac{1}{|I|}\int_I \omega\right)\left(\frac{1}{|I|}\int_I \omega^{-1/(p-1)}\right)^{p-1} \leq C$ for all $I \in \mathcal{D}$.

A weight $\omega$ is in the *dyadic doubling $A_1$ class* $\mathbf{A}_1^d$ if $\omega$ is dyadic doubling, and there is a constant $C$ such that $\frac{1}{|I|}\int_I \omega \leq C\,\omega(x)$ for a.e. $x \in I$, for all $I \in \mathcal{D}$.

The canonical examples are $\omega(x) = |x|^\alpha$. In this case $\omega \in \mathbf{A}_\infty^d$ if and only if $\alpha > -1$, $\omega \in \mathbf{RH}_p^d$ if and only if $\alpha > -1/p$, for $1 < p < \infty$ $\omega \in \mathbf{A}_p^d$ if and only if $-1 < \alpha < p - 1$, and $w \in \mathbf{A}_1^d$ if and only if $-1 < \alpha \leq 0$.

The class $\mathbf{A}_\infty^d$ is the union of the $\mathbf{RH}_p^d$ classes, and also of the $\mathbf{A}_p^d$ classes. The class $\mathbf{A}_1^d$ is strictly contained in the intersection of the $\mathbf{A}_p^d$ classes. More precisely:

$$(2.1) \qquad \mathbf{A}_\infty^d = \bigcup_{p>1} \mathbf{RH}_p^d = \bigcup_{p>1} \mathbf{A}_p^d, \quad \mathbf{A}_1^d \subset \bigcap_{p>1} \mathbf{A}_p^d.$$

These properties are known for weights that satisfy the conditions above for all intervals $I$; see [GC-RF]. In our case we consider only dyadic intervals, but we also assume that the weights are dyadic doubling, which yields (2.1); see [B]. An explicit example of a function in $\bigcap_{p>1} \mathbf{A}_p^d$ but not in $\mathbf{A}_1^d$ can be found in [JN].

To each weight $\omega \in \mathbf{A}_\infty^d$, with mean value one on $J_0$, we associate a function $b = b_\omega$

so that

(2.2) $$\omega(x) = \prod_{I \in \mathcal{D}} \Big(1 + b_I h_I(x)\Big), \qquad b_I = \frac{\langle \omega, h_I \rangle}{m_I \omega},$$

and
$$b(x) = \sum_{I \in \mathcal{D}} b_I h_I(x).$$

**Theorem 2.1 (R. Fefferman, Kenig, Pipher).** *Let $\{b_I\}_{I \in \mathcal{D}(J_0)}$ and $\varepsilon > 0$ be given, with $|b_I h_I| \leq (1-\varepsilon)$ for all $I \in \mathcal{D}$. Then the product (2.2) belongs to $\mathbf{A}_\infty^d$ if and only if $b \in \mathbf{BMO}^d$, that is, if and only if there is a constant $C$ such that*

(2.3) $$\sum_{I \in \mathcal{D}(J)} b_I^2 \leq C\,|J|, \quad \forall\, J \in \mathcal{D}.$$

*If so, the weight $\omega$ is a dyadic doubling weight.*

This is proved in [FKP].

The dyadic $\mathbf{RH}_p^d$ and $\mathbf{A}_p^d$ classes can also be characterized by summation conditions. Part (a) of the next theorem, and related results, appear in [B].

**Theorem 2.2.** *Let $\omega$ be a dyadic doubling weight $\omega$.*

(a) *(Buckley) $\omega \in \mathbf{A}_p^d$, $1 < p < \infty$, if and only if there is a constant $C$ such that*

$$\sum_{I \in \mathcal{D}(J)} \left(\frac{m_I \omega}{m_J \omega}\right)^{\frac{-1}{p-1}} b_I^2 \leq C\,|J|, \qquad \forall\, J \in \mathcal{D}.$$

(b) *$w \in \mathbf{A}_1^d$ if and only if there is a constant $C$ such that*

$$\frac{m_J \omega}{m_I \omega} \leq C, \quad \forall\, I \in \mathcal{D}(J),\ \forall\, J \in \mathcal{D}.$$

**Proof of (b):** ($\Rightarrow$) Integrate over $I \in \mathcal{D}(J)$ in the definition of $\mathbf{A}_1^d$. ($\Leftarrow$) By the Lebesgue differentiation theorem, the limit of $m_I \omega$ as a sequence of intervals $I \in \mathcal{D}(J)$ shrinks to a point $x \in J$, is $\omega(x)$ for almost every $x$.  □

Let $\widetilde{I}$ be the parent of $I$. Denote by $s_I$ the proportion of the mass of $\widetilde{I}$ that is carried by $I$. Then

(2.4) $$2 s_I = \frac{m_I \omega}{m_{\widetilde{I}} \omega} = 1 + b_{\widetilde{I}} h_{\widetilde{I}}(x), \qquad \forall\, x \in I.$$

**2.2. The $\lambda$-operation.** The nonlinear operation that sends the weight $\omega$ given by $\omega = \prod_{I \in \mathcal{D}}(1 + b_I h_I)$ into the weight $\omega_\lambda$,

$$\omega \mapsto w_\lambda = \prod_{I \in \mathcal{D}}(1 + \lambda b_I h_I), \tag{2.5}$$

is called the $\lambda$-*operation*.

*Remark* 2.3. Dyadic doubling and $\mathbf{A}_\infty^d$ weights are preserved under the $\lambda$-operation for $-1 \leq \lambda \leq 1$. This is a consequence of Theorem 2.1, and the observation that since $|\lambda b_I| \leq |b_I|$, $b \in \mathbf{BMO}^d$ implies $\lambda b \in \mathbf{BMO}^d$; see [P2].

To understand how the $\lambda$-operation affects $\mathbf{RH}_p^d$ and $\mathbf{A}_p^d$ weights, let us first consider how it affects quotients of mean values of $\omega$ over consecutive nested intervals. Accordingly define

$$s_I(\lambda) = \frac{m_I \omega_\lambda}{2\, m_{\widetilde{I}} \omega_\lambda}. \tag{2.6}$$

The $\lambda$-operation is a nonlinear operation on the weight, but it is linear at the level of the $s_I$'s:

**Lemma 2.4.** $s_I(\lambda) = \frac{1}{2} + \lambda(s_I - \frac{1}{2})$.

**Proof:** This is an immediate consequence of (2.4). □

**Lemma 2.5.** *Let $a_1, \ldots, a_n$ be positive numbers. Let $a_j(\lambda) = 1 + \lambda(a_j - 1)$, and $A_n(\lambda) = \prod_{j=1}^n a_j(\lambda)$. Then*

$$A_n(\lambda) \geq \min\{1, A_n(1)\}. \tag{2.7}$$

**Proof:** On expanding the product we get

$$A_n(\lambda) = \prod_{k=1}^n \left((1-\lambda) + \lambda a_k\right) = \sum_{i=0}^n (1-\lambda)^{n-i} \lambda^i \sum_{k_1 \neq \ldots \neq k_i} \prod_{j=1}^i a_{k_j}, \tag{2.8}$$

where the index notation $k_1 \neq \ldots \neq k_i$ means that the $k_j$'s are pairwise different, and $1 \leq k_j \leq n$. Now, using repeatedly the fact that for positive numbers the arithmetic mean is at least the geometric mean, we find that

$$\sum_{k_1 \neq \ldots \neq k_i} \prod_{j=1}^i a_{k_j} \geq N_{n,i} \left(\prod_{k_1 \neq \ldots \neq k_i} \prod_{j=1}^i a_{k_j}\right)^{1/N_{n,i}}, \tag{2.9}$$

where $N_{n,i}$ is the number of ways one can choose $i$ numbers from a collection of $n$ numbers, that is, $N_{n,i} = \frac{n!}{i!(n-i)!}$. Observe that in the double product in (2.9), each $a_j$ occurs as many times as possible ways one can choose $i-1$ numbers from a collection of $n-1$ numbers, that is $N_{n-1,i-1}$. Since $N_{n-1,i-1}/N_{n,i} = i/n$, the right hand side of (2.9) is equal to:

$$(2.10) \qquad N_{n,i}\left(\prod_{j=1}^n a_j\right)^{i/n} = N_{n,i}\left(A_n(1)\right)^{i/n}.$$

We conclude that

$$A_n(\lambda) \geq \sum_{i=0}^n (1-\lambda)^{n-i}\lambda^i N_{n,i}\left(A_n(1)\right)^{i/n} = \left(\lambda\left(A_n(1)\right)^{1/n} + (1-\lambda)\right)^n.$$

Fix $A \geq 0$. The function $f(\lambda) = \left(\lambda A + (1-\lambda)\right)^n$ is monotonic for $\lambda \in [0,1]$, hence it is bounded below by $\min\{f(0), f(1)\}$. Set $A = \left(A_n(1)\right)^{1/n}$; then $f(0) = 1$, $f(1) = A_n(1)$. This completes the proof of the lemma. $\square$

## 3. The theorem for $\mathbf{A}_p^d$ weights

**Theorem 3.1.** *Given a weight $\omega \in \mathbf{A}_p^d$, then $\omega_\lambda \in \mathbf{A}_p^d$ for all $1 \leq p \leq \infty$ and all $0 \leq \lambda \leq 1$.*

**Proof:** By Remark 2.3, $\omega_\lambda$ is a dyadic doubling weight, and $w \in \mathbf{A}_\infty^d$ implies $\omega_\lambda \in \mathbf{A}_\infty^d$.

For dyadic intervals $I \subset J$, the ratios $\frac{m_I\omega}{m_J\omega}$ and $\frac{m_I\omega_\lambda}{m_J\omega_\lambda}$ can be explicitly computed as products of $s_K$'s and $s_K(\lambda)$'s respectively:

$$(3.1) \qquad \frac{m_I\omega}{m_J\omega} = \prod_{I \subseteq K \subset J} 2s_K, \quad \text{and} \quad \frac{m_I\omega_\lambda}{m_J\omega_\lambda} = \prod_{I \subseteq K \subset J} 2s_K(\lambda),$$

where $K \in \mathcal{D}(J)$, $2s_K(\lambda) = 1 + \lambda(2s_K - 1)$, and $2s_K = m_K\omega/m_{\widetilde{K}}\omega$.

These are the kind of products considered in Lemma 2.5. Since $\omega$ is dyadic doubling, (2.4) shows that $\varepsilon \leq 2s_K \leq 2 - \varepsilon$. In particular $a_K = 2s_K > 0$. We can apply Lemma 2.5 to each pair $(I,J)$ of dyadic intervals, $I \subset J$, where $n = \log_2(|J|/|I|)$, $A_n(\lambda) = m_I\omega_\lambda/m_J\omega_\lambda$, and $A_n(1) = m_I\omega/m_J\omega$, to conclude that

$$(3.2) \qquad \left(\frac{m_I\omega_\lambda}{m_J\omega_\lambda}\right)^r \leq \max\left\{1, \left(\frac{m_I\omega}{m_J\omega}\right)^r\right\},$$

for all negative $r$.

**Case $1 < p < \infty$:** By Theorem 2.2(a), it suffices to find a constant $C$ such that

$$(3.3) \qquad \frac{1}{|J|}\sum_{I \in \mathcal{D}(J)} \left(\frac{m_I\omega_\lambda}{m_J\omega_\lambda}\right)^{\frac{-1}{p-1}} b_{I,\lambda}^2 \leq C \quad \forall J \in \mathcal{D},$$

where $b_{I,\lambda} = \langle \lambda b, h_I \rangle = \lambda b_I$.

We use (3.2), with $r = -1/(p-1)$, to estimate the left hand side of (3.3). Using the observation that $b_{I,\lambda}^2 \leq b_I^2$ for $0 < \lambda < 1$, we get

$$\frac{1}{|J|} \sum_{I \in \mathcal{D}(J)} \left(\frac{m_I \omega_\lambda}{m_J \omega_\lambda}\right)^{\frac{-1}{p-1}} b_{I,\lambda}^2 \leq \frac{1}{|J|} \sum_{I \in \mathcal{D}(J)} \max\left\{1, \left(\frac{m_I \omega}{m_J \omega}\right)^{\frac{-1}{p-1}}\right\} b_I^2$$

$$\leq \frac{1}{|J|} \sum_{I \in \mathcal{D}(J)} b_I^2 + \frac{1}{|J|} \sum_{I \in \mathcal{D}(J)} \left(\frac{m_I \omega}{m_J \omega}\right)^{\frac{-1}{p-1}} b_I^2 \leq C;$$

where $C$ is a constant independent of the interval $J$. The last inequality holds since $\omega \in \mathbf{A}_p^d \subset \mathbf{A}_\infty^d$, by Theorems 2.1 and 2.2(a).

**Case $p = 1$:** By Theorem 2.2(b), it suffices to find a constant $C$ such that

$$(3.4) \qquad \left(\frac{m_I \omega_\lambda}{m_J \omega_\lambda}\right)^{-1} \leq C, \quad \forall I \in \mathcal{D}(J), \quad \forall J \in \mathcal{D}.$$

We use (3.2) again, with $r = -1$, to conclude that the left hand side of (3.4) is bounded by $\max\left\{1, \left(\frac{m_I w}{m_J w}\right)^{-1}\right\}$, which in turn is bounded because $\omega \in \mathbf{A}_1^d$. □

## 4. Examples in $\mathbf{RH}_p^d$

**Theorem 4.1.** *For each $p > 1$, there exist a dyadic doubling weight $\omega$ in $\mathbf{RH}_p^d$ and a number $\lambda \in (0,1)$ such that $\omega_\lambda$ is not in $\mathbf{RH}_p^d$.*

**Proof:** Our examples are of the following form. Let $I_i = [2^{-i}, 2^{-i+1}]$ for $i \geq 1$, and $J_i = [0, 2^{-i}]$ for $i \geq 0$. Fix numbers $s_i \in (0,1)$ for $i \geq 1$. For $x \in [0,1]$ let

$$(4.1) \qquad \omega(x) = \sum_{i=1}^{\infty} c_i \chi_{I_i}(x), \quad \text{where} \quad c_i = 2^i \, s_1 \ldots s_{i-1}(1 - s_i).$$

One can think of $\omega$ as the weight of mass one obtained by assigning the fractions $s_i$ and $1 - s_i$ of the mass of $J_{i-1}$ to the left half $J_i$ and the right half $I_i$ respectively of $J_{i-1}$, for each $i \geq 1$. In the notation of Section 2,

$$(4.2) \qquad s_i = s_{J_i} = \frac{m_{J_i} \omega}{2 \, m_{J_{i-1}} \omega}, \quad \text{and} \quad m_{J_i} \omega = 2^i \, s_1 \ldots s_{i-1} s_i.$$

We further assume that the sequence $\{s_i\}$ is *n-periodic*: there is a smallest positive integer $n$ such that $s_{i+n} = s_i$ for all $i \geq 1$. Such a $\omega$ is a dyadic doubling weight. For $n$-periodic weights of this form,

$$(4.3) \qquad c_{kn+j} = (2^n \, s_1 \ldots s_n)^k \, c_j, \qquad \forall \, k \geq 0, j \geq 1.$$

**Claim:** *The* $\mathbf{RH}_p^d$ *condition* $\frac{1}{|I|}\int_I \omega^p \leq C\left(\frac{1}{|I|}\int_I \omega\right)^p$, *for all* $I \in \mathcal{D}$, *for an n-periodic weight of the form (4.1) reduces to the condition*

$$(4.4) \qquad 2^n s_1 \ldots s_n < 2^{n/p}.$$

To see this, first note that $\omega$ is constant on each dyadic interval $I$ whose left endpoint is not 0, so the reverse Hölder $p$ condition holds with $C = 1$ on these intervals. The only other dyadic intervals are $J_l = [0, 2^{-l}]$, $l \geq 0$.

Write $l = qn + m$ where $q$ and $m$ are non-negative integers and $0 \leq m \leq n-1$. We use (4.3) to compute the mean of $\omega^p$ on $J_l = \bigcup_{i=l+1}^\infty I_i$:

$$\begin{aligned}
\frac{1}{|J_l|}\int_{J_l} \omega^p &= 2^l \sum_{i=l+1}^\infty \int_{I_i} \omega^p = 2^l \sum_{i=qn+m+1}^\infty c_i^p \, 2^{-i} \\
&= 2^l \sum_{k=q}^\infty \sum_{j=1}^n (c_{kn+j+m})^p \, 2^{-(kn+j+m)} \\
&= 2^l \sum_{k=q}^\infty \sum_{j=1}^n (2^n s_1 \ldots s_n)^{kp} (c_{j+m})^p \, 2^{-(kn+j+m)} \\
(4.5) \qquad &= 2^l \, d_{p,m} \sum_{k=q}^\infty \left[(2^n s_1 \ldots s_n)^p \, 2^{-n}\right]^k,
\end{aligned}$$

where

$$d_{p,m} = \sum_{j=1}^n (c_{j+m})^p \, 2^{-(j+m)}.$$

The series (4.5) converges if and only if $2^n s_1 \ldots s_n < 2^{n/p}$; this is condition (4.4).

The mean of $\omega$ on $J_l$ is

$$(4.6) \qquad \frac{1}{|J_l|}\int_{J_l} \omega = 2^l \, d_{1,m} \sum_{k=q}^\infty (s_1 \ldots s_n)^k;$$

this series converges since each $s_i$ is less than one.

If $2^n s_1 \ldots s_n < 2^{n/p}$, then on summing the series in (4.5) and (4.6) we find

$$(4.7) \qquad \frac{\frac{1}{|J_l|}\int_{J_l} \omega^p}{\left(\frac{1}{|J_l|}\int_{J_l} \omega\right)^p} = 2^{l(1-p)} \frac{d_{p,m}}{(d_{1,m})^p} \, 2^{n(p-1)q} \, \frac{[1-(s_1 \ldots s_n)]^p}{1-(2^n s_1 \ldots s_n)^p \, 2^{-n}},$$

for all $l = qn + m \geq 0$ with $q \geq 0$ and $0 \leq m \leq n-1$. This expression depends on $p$, $n$, and $m$ but is independent of $q$, since $d_{p,m}$ and $d_{1,m}$ are independent of $q$ and the exponent of 2 in the right hand side of (4.7) reduces to $m(1-p)$. Taking the maximum over $m \in \{0, \ldots, n-1\}$, we see that there is a uniform upper bound, depending only on $p$ and $n$, for $\frac{1}{|J_l|}\int_{J_l} \omega^p / (\frac{1}{|J_l|}\int_{J_l} \omega)^p$, for all $l \geq 0$. Therefore $\omega \in \mathbf{RH}_p^d$, which establishes the claim.

Let $f(x_1, \ldots, x_n) = 2^n x_1 \ldots x_n$. Then by the claim above, the $n$-periodic weight $\omega$ associated to the point $(s_1, \ldots, s_n)$ in the open unit $n$-cube $(0,1)^n$ is in $\mathbf{RH}_p^d$ if and only if $f(s_1, \ldots, s_n) < 2^{n/p}$.

By Lemma 2.4, the $\lambda$-operation takes the weight $\omega$ with sequence $\{s_i\}$ to a new weight $\omega_\lambda$ with sequence $\{s_i(\lambda)\}$, where

$$s_i(\lambda) = \frac{1}{2} + \lambda \left( s_i - \frac{1}{2} \right). \tag{4.8}$$

In the unit $n$-cube, the $\lambda$-operation moves the point $P : (s_1, \ldots, s_n)$ towards the centre $Q : (1/2, \ldots, 1/2)$ of the cube along the line segment joining $P$ to $Q$. For $\lambda \in (0,1)$, the point $P_\lambda : \big(s_1(\lambda), \ldots, s_n(\lambda)\big)$ associated to $\omega_\lambda$ lies on $L$ between $P$ and $Q$.

We wish to find a weight $\omega$ and a number $\lambda \in (0,1)$ such that $\omega \in \mathbf{RH}_p^d$ but $\omega_\lambda \notin \mathbf{RH}_p^d$. To do this, we will exhibit a line $L$ through $Q$, and points $P$ and $P_\lambda$ on $L$, such that: $(i)$ $P$ and $P_\lambda$ are inside the open unit cube $(0,1)^n$, and $P_\lambda$ is between $P$ and $Q$; $(ii)$ the value of $f$ at $P$ is strictly less than $2^{n/p}$; and $(iii)$ the value of $f$ at $P_\lambda$ is at least $2^{n/p}$. Then: $(i)'$ there is a unique $\lambda \in (0,1)$ such that the $\lambda$-operation takes the weight $\omega$ associated to $P$ to the weight $\omega_\lambda$ associated to $P_\lambda$; $(ii)'$ the weight $\omega$ associated to $P$ is in $\mathbf{RH}_p^d$; and $(iii)'$ the weight $\omega_\lambda$ associated to $P_\lambda$ is not in $\mathbf{RH}_p^d$.

Fix $a \in (0, 1/2)$. Let $L$ be the line through $(a, 1, \ldots, 1)$ and $Q$; it has equation $\left(a + \left(\frac{1}{2} - a\right)t, 1 - \frac{t}{2}, \ldots, 1 - \frac{t}{2}\right)$, and values of $t \in [0, 2]$ correspond to points in the closed unit cube $[0,1]^n$. Let $g(t)$ be the restriction of $f$ to $L$, for $t \in [0, 2]$:

$$g(t) = \big(2a + (1 - 2a)t\big)(2 - t)^{n-1}. \tag{4.9}$$

Then $g'(t) = (2 - t)^{n-2} \left[ -n(1 - 2a)t + 2(1 - (n+1)a) \right]$.

If $n > 2$, $g$ has two critical points; if $n = 2$, $g$ has a single critical point. Since $g(2) = f(1 - 2a, 0, \ldots, 0) = 0$, and $f$ is positive in the open unit cube, a critical point of $g$ at $t = 2$ is a minimum. The critical point which occurs at

$$t_m = \frac{2(1 - (n+1)a)}{n(1 - 2a)} \tag{4.10}$$

is a maximum for $g$. Denote this point by $P_m$. The maximum value is

$$g(t_m) = \frac{2^n (n-1)^{n-1} (1 - a)^n}{n^n (1 - 2a)^{n-1}}. \tag{4.11}$$

This maximum is taken inside the open unit cube if and only if $t_m > 0$; in other words, if and only if $a < (n+1)^{-1}$. In this case, $P_m$ lies between $(a, 1, \ldots, 1)$ and the centre $Q$ of the cube.

The maximum value of $g$ is a continuous, increasing function of $a$, for $a \in [0, (n+1)^{-1}]$. We denote this function by $h(a) = g(t_m)$ (see (4.11) above). Then $h'(a) > 0$ for $a \in [0, (n+1)^{-1}]$. Also, $h(0) = 2^n (n-1)^{n-1} n^{-n}$, and $h\left(\frac{1}{n+1}\right) = \frac{2^n}{n+1}$. Therefore, if $p > 1$ is large enough that

$$2^{n/p} < \frac{2^n}{n+1}, \tag{4.12}$$

then we can choose $a_p \in (0, (n+1)^{-1})$ so that $h(a_p) \geq 2^{n/p}$ and so that the maximum value $h(a_p)$ of $f$ along the line segment from $(a_p, 1, \ldots, 1)$ to $Q$ is taken on at a point inside the open unit cube.

Notice that if $p$ is too large, we may not be able to find a line segment on which the maximum value $h(a_p)$ of $f$ is actually equal to $2^{n/p}$, since the smallest $h(a_p)$ is $h(0)$ which is strictly larger than one (if $n > 2$).

Define $P_\lambda$ as follows: if $h(0) < 2^{n/p} < \frac{2^n}{n+1}$, then choose the unique $a_p \in (0, (n+1)^{-1})$ such that $h(a_p) = 2^{n/p}$, and let $P_\lambda = P_m$ for that $a_p$. If $2^{n/p} \leq h(0)$, choose $a_p$ small enough that $f(a_p, 1, \ldots, 1) = 2^n a_p < 2^{n/p}$, and let $P_\lambda = P_m$ for that $a_p$.

In both cases, $f(a_p, 1, \ldots, 1) < 2^{n/p}$ while the value of $f$ at $P_m$ is at least $2^{n/p}$. Since $g$ is increasing for $0 < t < t_m$, we may choose a point $P$ between $(a_p, 1, \ldots, 1)$ and $P_\lambda$ such that the value of $f$ at $P$ is strictly less than $2^{n/p}$. Let $\lambda$ be the unique number in $(0, 1)$ such that the $\lambda$-operation takes this point $P$ to $P_\lambda$.

The weights $\omega$ and $\omega_\lambda$ associated to $P$ and $P_\lambda$ respectively, together with the number $\lambda$, furnish our example: $\omega$ is in $\mathbf{RH}_p^d$ but $\omega_\lambda$ is not in $\mathbf{RH}_p^d$.

The critical value of $p$, from (4.12), is

$$p_c = \frac{n \log 2}{n \log 2 - \log(n+1)}. \tag{4.13}$$

Examples for all $p > p_c$ can be found among the $n$-periodic weights of the form described above. Letting $n$ tend to infinity, we see that $p_c$ tends to 1, and so we have found examples for all $p > 1$. $\square$

CRISTINA PEREYRA, DEPARTMENT OF MATHEMATICS, PRINCETON UNIVERSITY, PRINCETON NJ 08544-0001
  *E-mail address*: crisp@math.princeton.edu

LESLEY WARD, DEPARTMENT OF MATHEMATICS, RICE UNIVERSITY, HOUSTON TX 77251-1892; AND MSRI, 1000 CENTENNIAL DRIVE, BERKELEY CA 94720-5070
  *E-mail address*: lesley@math.rice.edu